\documentclass{article}
\newcommand{\R}{I\!\! R}
\newcommand{\C}{I\!\! C}

\newcommand{\N}{I\!\! N}

\begin{document}

Tsemo, Aristide

Department of Mathematics and Computer Sciences

Ryerson University

350, Victoria Street M1L 2V5, Toronto Canada

tsemo58@yahoo.ca

\bigskip
\bigskip
\bigskip
\centerline{\bf Deformation of Homogeneous structures and
Homotopy of symplectomorphisms groups.}

\bigskip
\bigskip

{\bf Introduction.}

\bigskip

Let $(N,\omega)$ a symplectic manifold, the group $Symp(N)$ of
symplectomorphisms of $(N,\omega)$ acts transitively on it.
Moreover, Banyaga has shown that $Sym(N)$ determines completely the
symplectic structure of $(N,\omega)$   when $N$ is compact. This
motivates the study of the properties of $Symp(N)$ which must enable
to understand the geometric properties of $(N,\omega)$.

 The beginning
of the study of the homotopy properties of $Symp(N)$ has its origin
in the theory of pseudo-holomorphic curves defined by Gromov, using
this theory, Gromov has shown that the group of compactly supported
symplectomorphisms in the interior of the symplectic $4$-dimensional
standard ball is contractible. Let $(S^2\times S^2,\omega_0)$ be the
product of two copies of the $2$-dimensional sphere $S^2$, endowed
with the symplectic structure which is the product of the canonical
volume form of $S^2$. Gromov has also shown that $Symp(S^2\times
S^2)$ has the same homotopy type than $SO(3)\times SO(3)$.

\bigskip

 The work of Gromov has been generalized by many authors,
remark that if we endow $S^2\times S^2$ with a symplectic form
$\omega$ such that  the volume of the fiber of the fibration
$S^2\times S^2\rightarrow S^2$ is not the volume of its base, the
homotopy of $Symp(S^2\times S^2,\omega)$ can be different than the
homotopy of a finite dimensional Lie group. This phenomenon has been
firstly observed by Abreu, and completely studied by McDuff and
Abreu. Their study use the theory of pseudo-holomorphic curves, and
stratification of  the space of pseudo-convex structures.

In his thesis Pinsonnault has studied the homotopy type of the one
point blow-up $N$ of $S^2\times S^2$, if the volume of the fiber and
base of the fibration $S^2\times S^2$ are equal, he has shown that
$Symp(N)$ has the homotopy type of a $2$-dimensional torus. Joseph
Coffey has used the decomposition of a symplectic $4$-dimensional
manifold to relate the homotopy type of the group of
symplectomorphisms $Symp(N)$ of a $4$-dimensional manifold to a
configuration space, more precisely Paul Biran has shown that we can
remove a simplex $D$ on a $4$-dimensional symplectic manifold in
such a way that the resulting space is a disc bundle over a
symplectic surface, the group $Symp(N)$ is then an extension of its
subgroup which preserves $D$ by the orbit of $D$ by the group
$Symp(N)$. This approach must be related to the work of Lalonde and
Pinsonnault who have studied the relations of the space of embedding
of a symplectic ball with a given radius in a symplectic
$4$-dimensional manifold and $Symp(N)$.

\bigskip
 As remarked Pinsonnault in his thesis, the calculation of
the homotopy type of the group of symplectomorphisms of $S^2\times
S^2$ and its blow-up is possible because the topology of these
manifolds is very simple, and there exists a good knowledge of the
theory of pseudo-holomorphic curves of these manifolds. Another
simple example of symplectic $4$-dimensional manifold is the
$4$-dimensional torus $N^4$, the study of $Symp(N^4)$ is hardly
tractable using the theory of pseudo-holomorphic curves since
$H^2(N^4,{\R})$ is a $6$-dimensional real vector space.

\medskip

 The purpose of this paper is to describe the homotopy type
of $Symp(N^4)$ using new ideas. The torus $N^4$ is the quotient of
the affine space ${\R}^4$ by $4$-translations which preserve the
standard symplectic form. Let $D'$ be the boundary of a fundamental
domain of this action. It projects to $N^4$ to define a union $D$ of
$3$-dimensional torus. Our approach may be related to the work of
Joseph Coffey and Lalonde-Pinsonnault, we remark that the group
$Symp(N^4)$ is the total space of a fibration whose base space is
the orbit of $D$ under $Symp(N^4)$ and fiber the group of
symplectomorphisms which preserve $D$.

  To determine the homotopy type of $Symp(N^4)$, one needs to
determine the homotopy type of the space $N(D)$ of orbits of $D$
under $Symp(N)$. This is performed using the theory of deformation
of homogeneous structures as it is described by W.Goldman. We show
that
 $N(D)$ is  contractible.

\bigskip
{\bf Decomposition of symplectic manifolds and group of
symplectomorphism.}

\medskip
The purpose of this paragraph is to define a notion of canonical
decomposition of a compact symplectic $4$-dimensional manifold, and
to show that the connected component of the group of Hamiltonian
diffeomorphisms which preserve such a configuration is contractible.

\medskip
 {\bf Definition}

Let $(N,\omega)$ be a compact $4$-dimensional symplectic manifold, a
configuration $D$ of $(N,\omega)$ is a finite union $D$ of compact
$3$-dimensional submanifolds $D_1,..,D_p$ which satisfy the
following conditions:

The manifold $N-D$ is symplectomorphic to the $4$-dimensional
symplectic standard ball or to a polydisc endowed with the
restriction of the standard symplectic form of ${\R}^4$, $D$ is
connected.

The restriction of $\omega$ defines on $D_i$ a contact structure
such that the fibers of the characteristic foliation are circles,
$D_i$ is the total space of a bundle whose leaves are the
characteristic leaves, and which base space is an intersection
$D_i\cap D_j$. The fibers of the characteristic foliation of $D_i$
intersects at least two different cells.

The intersection $D_{ij}=D_i\cap D_j$ is either a symplectic
$2$-dimensional submanifold of $(N,\omega)$, or a Lagrangian
$2$-dimensional submanifold stable by the characteristic foliations.

Suppose that $D_{ij}$ is a symplectic submanifold and consider $h_t$
a path of Hamiltonian diffeomorphisms of  $D_{ij}$ such that the
restriction of $h_t$ to $D_{i_1i_2i_3}=D_{i_1}\cap D_{i}\cap D_{j}$
is the identity if $D_{i_1i_2i_3}$ is a circle, then we can extend
$h_t$ to a path of Hamiltonian diffeomorphisms $h'_t$ of $N$ in such
a way that it preserves the configuration, and its restriction to
the cells $D_l$ such that $D_l\cap D_{ij}$ is either a circle, or is
empty is the identity. The restriction of $h'_t$ to a symplectic
surface $D_{i'j'}$ different of $D_{ij}$ is the identity.

 The configuration defines on $N$ the structure of a CW-complex
with one $4$-dimensional cell.

An example of a configuration is the $4$-dimensional torus $N^4$. It
is the product of $2$-dimensional torus $N_2$ and $N'_2$. $N^4$ is
endowed with the symplectic form which is the product of the
standard symplectic forms of $N_1$ and $N_2$ Let $l_1$, $l_2$, be
the two curves of $N_2$ parallel in respect to the flat structure
which generate $\pi_1(N_2)$, and $l_3$, $l_4$, the curves of $N'_2$
which satisfy the same properties. The configuration is the union of
$N_2\times l_3, N_2\times l_4, N'_2\times l_1, N'_2\times l_2$. Each
function $H:N_2\rightarrow {\R}$, which is constant on $l_1$ and
$l_2$ can be extended to $N_2\times N'_2$ by the projection
$N_2\times N'_2\rightarrow N_2$, and it preserves the configuration.

\medskip

We denote by $Symp(N,D)$ the group of symplectomorphisms of
$(N,\omega)$ which preserve $D$,  $Ham(N,D)$ its subgroup of
Hamiltonian symplectomorphisms, $Symp_{Id}(N,D)$ is the subgroup of
$Symp(N,D)$ whose restriction to $D$ is the identity, and
$Ham_{Id}(N,D)$ is the subgroup of Hamiltonian diffeomorphisms of
$Symp_{Id}(N,D)$.

\medskip

{\bf Proposition} {\it The connected component $Ham_{Id}(N,D)_0$ of
the group $Ham_{Id}(N,D)$ is  contractible.}

\medskip

 A theorem of Weinstein asserts that there exists a neighborhood
$U_{Id}$ of the identity in $Symp(N)$ which is contractible. Without
restricting the generality, we can suppose that the image of
$U_{Id}$ by the flux map is a contractible open subset. Thus we can
suppose that $Ham(N)\cap U_{Id}$ is contractible, and the existence
of a continuous map $Ham(N)\cap U_{Id}$ to the space of exact time
dependent $1$-forms, $h\rightarrow dH^h$, where $H^h:N\rightarrow
{\R}$ is a time dependent differentiable function such that $h$ is
the value at $1$ of the flow generated by $dH$.

Without restricting the generality, we can suppose that for every
element $h\in Ham_{Id}(N,D)_0$, there exists a  neighborhood  $V$ of
$D$ such that the restriction of $h$ to $V$ is the restriction of an
element of $U_{Id}$, $h'$ to $V$ which coincide with $h$ on an open
neighborhood $U$ which contains $V$, and $h(V)\subset U_1\subset U$,
where $U_1$ is open. The existence of $h'$ can be shown using a cut
function. We can also suppose that and $N-V,N-U$ are polydiscs. Let
$f$ be a cut function such that the restriction of $f$ to $U_1$ is
$1$ and the restriction of $f$ to $N-U$ is zero, we define a map
$\Phi_t$, $h\rightarrow {\psi_t}^{-1}\circ h$, where $\psi_t$ is the
flow generated by $fH^{h'}$, $\Phi_0(h)=h$, $\Phi_1(h)$ is an
Hamiltonian map whose support is contained in the interior of $N-V$.
Let $Ham_{Id}(N,V,U_1,U,D)$ be the subset of $Ham_{Id}(N,D)_0$ such
that the restriction of $h$ to $V$   is the restriction of an
element $h'$ of $U_{Id}$, and $h(V)\subset U_1$, the restrictions of
$h$ and $h'$ to $U$ coincide. A result of Gromov says that the group
of symplectomorphisms whose supports are in the interior of $N-V$ is
contractible. The map $\Phi_t$ implies that $Ham_{Id}(N,V,U_1,U,D)$
is contractible.

Let  ${(V_n)}_{n\in {\N}}$, be a family of open subsets which
contain $D$ such that $V_{n+1}$ is contained in $V_n$,
$\cap_{n\in{\N}} V_n=D$. Consider the families of open subsets $U_n$
and $U^n_1$, such that $V_n\subset U^n_1\subset U_n$, we suppose
that $\cap_{n\in{\N}}U_n=\cap_{n\in{\N}}U^n_1=D$,
$U^{n+1}_{1}\subset U^n_1$, $U_{n+1}\subset U_n$. Then
$Ham_{Id}(N,D)_0$ is the  limit of
$Ham_{Id}(N,V_{n_1},U^{n_2}_1,U_{n_3},D)$ when the previous
expression has a sense. We suppose that $N-V_n$ is a polydisc. Since
$V_{n+1}$ is contained in the interior of $V_n$, we deduce that
$Ham_{Id}(N,V_{n},U^n_1,U_n,D)$ is contained in the interior of
$Ham_{Id}(N,V_{n+1},U^{n+1}_1,U_n,D)$. This implies that
$Ham_{Id}(N,D)_0$ is contractible.

\medskip

{\bf Proposition.} {\it The group  $Ham(N,D)_0$ is  contractible.}

\medskip

{\bf Proof.} Let $L_D$ be the restriction to $D$ of the subgroup
$Ham'(N,D)$ of $Ham(N,D)_0$ whose restriction to $D_{i_1i_2i_3}$ is
the identity if $D_{i_1i_2i_3}$ is a curve. If $D_{ij}$ is a
symplectic submanifold, we denote by $L_{ij}$ the group of
Hamiltomian diffeomorphisms of $D_{ij}$ whose restriction to
$D_{i_1i_2i_3}$ is the identity, if $D_{i_1i_2i_3}$ is a cercle, we
denote by $L'_D$ the product of the groups $L'_{ij}$. The extension
property verified by the configuration implies that the restriction
map $L_D\rightarrow L'_D$ is a fibration. The fiber of this
fibration is contractible, since the gauge group of the fibration
induced by the characteristic foliation is contractible, since each
fibers intersects two different cells, and a connected subgroup of
diffeomorphisms of an interval is contractible. Since $L'_D$ is
contractible, we deduce that $L_D$ is contractible.

Let $L^1_D$ be the restriction of $Ham(N,D)_0$ to $D$, and $L^2_D$
the restriction of $Ham(N,D)_0$ to the circles $D_{i_1i_2i_3}$, we
have a fibration $L^1_D\rightarrow L^2_D$ whose fiber is $L_D$,
since $L_D$ is contractible, and $L^2_D$ is contractible since a
connected group of diffeomorphism of an interval is contractible, we
deduce that $L^1_D$ is a contractible. The fiber of the fibration
$Ham(N,D)_0\rightarrow L^1_D$ is $Ham_{Id}(N,D)_0$ which is
contractible. This implies that $Ham(N,D)_0$ is contractible.

\bigskip

{\bf Isotopy of configurations and homogeneous structures.}

\medskip
To determine the homotopy of $Ham(N)$, we have to determine the
homotopy of the orbit of $D$ under the action of $Ham(N)$. To solve
this problem, we consider the situation when $N$ is a locally
homogeneous manifold.
\medskip

{\bf Definition} Let $V$ be a differentiable manifold, and $H$ a Lie
subgroup acting on $V$, the action of $H$ satisfies the unique
extension property, if and only if two elements $h_1$, $h_2$ of $H$
which coincide on an open subset of $V$,  coincide on $V$.
 \medskip

{\bf Definition}
 Let $(V,H)$ be a $n$-dimensional  manifold endowed with the action of a Lie
 group which satisfies the unique extension property. A manifold $N$ is endowed
 with a $(V,H)$ (homogeneous) structure, if and only if there exists
 a $V$-atlas $(U_i,\phi_i)_{i\in I}$ of $N$, such that $\phi_i:U_i\rightarrow V$,
 and $\phi_j\circ{\phi_i}^{-1}$ is the restriction of the action
 of an element $h_{ij}$ of $H$ on $\phi_i(U_i\cap U_j)$.
 \medskip

 Let $N$ be a manifold endowed with a $(V,H)$ homogeneous structure,
 and $\hat N$ be the universal cover $N$, $\hat N$ inherits from $N$
 a structure (the pull-back) of a $(V,H)$ manifold defined by a local diffeomorphism
 $D_N:\hat N\rightarrow V$ called the developing map. This
 developing map gives rise to a representation
 $h_N:\pi_1(N)\rightarrow H$ defined by the following commutative
 diagram:

 \begin{eqnarray*}
 \matrix{\hat N &{\buildrel{d}\over{\longrightarrow}}&\hat N\cr
 \downarrow D_N&&\downarrow D_N\cr
 V &{\buildrel{h_N(d)}\over{\longrightarrow}}&V}
 \end{eqnarray*}
called the holonomy representation.

\medskip

Examples of homogeneous manifolds are:

Affine manifolds, here $V$ is the affine space ${\R}^n$, and $H$ is
$Aff({\R}^n)$ the group of affine transformations.

Projective manifolds, here $V$ is either $P^n({\R})$, or
$P^n({\C})$, and $H$ is the group of projective transformations.

\bigskip

Let $N$ be a manifold endowed with a $(V,H)$ structure, we denote by
$D(N,V,H)$ the space of $(H,N)$ structures defined on $N$. To each
element $[h_N]$ of $D(N,V,H)$, we can associate the following flat
bundle $V_{h_N}$: let $(U_i,\phi_i)$ be the atlas which defines the
$(V,H)$-structure of $N$. The total space of $V_{h_N}$ is the
quotient of $\hat N\times V$ by the action of $\pi_1(N)$ defined on
the first factor by the Deck transformations, and on the second
factor by the holonomy $h_N$ of the $(V,H)$ structure. The following
result his shown in Goldman:

\medskip

{\bf Theorem.}

{\it Let $N$ be a compact manifold endowed with a $(V,H)$ structure,
then the map $D(N,V,H)\rightarrow Rep(\pi_1(N),H)$ which assigns to
each $(V,H)$ structure its holonomy representation is an open map,
two nearby $(V,H)$ structures with the same holonomy are isotopic.}

Suppose now that the $(V,H)$-manifold $N$ is endowed with a
symplectic form $\omega$, and there exists a configuration $D$ of
$N$ such that $N-D$ is symplectomorphic to the standard ball, or the
standard polydisc. The idea that will use is to study the
deformation space of $D$ under the action of the group of
symplectomorphisms, using the deformation theory of locally
homogeneous structures developed in Goldman. Let $\phi$ be a
symplectomorphism of $(N,\omega)$ the symplectomorphism $\phi$ acts
on $D(N,V,H)$ as follows: let $(U_i,\phi_i)_{i\in I}$ be an atlas
which defines the $(V,H)$ structure, then $(\phi(U_i),\phi_{i'}\circ
\phi)$ defines a $(V,H)$ structure $\phi([h_N])$. (without
restricting the generality, we suppose that there exists $i'\in I$
such that $\phi(U_i)$ is contained in the domain of a chart
$U_{i'})$.  The affine structure of $h_N$ is defined by the
projection on the trivialization $U_i\times V$, on the second
factor. It can be defined also by the following section $l_{h_N}$ of
$V_{h_N}$: we define the restriction $l^i_{h_N}$ of $l_{h_N}$ to
$U_i$ to be the map defined by $l^i_{h_N}(u)=(u,\phi_i(u))$.  The
affine structure of $\phi(h_N)$ is defined by a  section of
$V_{h_N}$ transverse to the horizontal foliation after the
identification between $V_{h_N}$ and $V_{\phi(h_N)}$. Recall that
the horizontal foliation of $V_{h_N}$ is the foliation whose leaves
are the projections of the submanifolds $\hat N\times \{u\}$ of
$\hat N\times V$ by the covering map $\hat N\times V\rightarrow
V_{h_N}$.

\medskip

Suppose that $(V,H)$ is $({\R}^n,Aff({\R}^n))$, in this situation,
the bundle $V_{h_N}$ is the affine vector bundle defined on $N$, it
is associated to a principal bundle defined on $N$, which is the
bundle of affine frames. Recall that an affine manifold is also
defined by a connection whose curvature and torsion forms vanish
identically, this connection induces on the bundle of affine frames
a connection a la Erhesmann, the horizontal distribution which
defines the Erhesmann connection induced the horizontal foliation of
$h_{h_N}$.

\bigskip

{\bf The action of  $\hat{Symp(N)}_0$ on $V_{h_N}$.}

\medskip

We are going to define an action of the connected component of the
universal cover $\hat{Symp(N)}_0$ of the group of symplectomorphisms
$Symp(N)$ of $(N,\omega)$ on $V_{h_N}$. We suppose that our $(V,H)$
structure is complete, or more generally, we suppose that the
developing map is a covering map, this is equivalent to saying that
$N$ is the quotient of an open subset $U_N$ of $V$ by a discrete
subgroup of $H$.

We define $V^1_{h_N}$ to be the $U_N$-bundle defined on $N$, by
making the quotient of $\hat N\times U_N$ by the action of
$\pi_1(N)$ defined on the first factor by the deck transformations,
and on the second factor by the action of the holonomy.

Let $h'$ be an element of ${Symp(N)}_0$, we can lift $h'$ to an
element $h"$ of $U_N$, since $h'$ is an element of $Symp(N)_0$, $h"$
commutes with the action of $\pi_1(N)$ on $U_N$. We denote by
$Symp'(N)_0$, the group of symplectomorphisms of $(U_N,\omega')$
(where $\omega'$ is the symplectic form of $U_N$ lifted to $U_N$ by
the covering map $p':U_N\rightarrow N$) which commutes with the
holonomy: this group is a covering of $Symp(N)_0$. Since the element
$h"$ commutes with $\pi_1(N)$, it induces a gauge transformation on
$V^1_{h_N}$ which covers the identity of $N$. The covering map
$\hat{Symp(N)}_0\rightarrow Symp'(N)$ induces an action of
$\hat{Symp(N)}_0$ on $V^1_{h_N}$.

To study the action of $Symp(N)$ on the configuration $D$, we make
the following assumption: let $N_U$ be a fundamental domain of the
action of $\pi_1(N)$ on $U_N$, we denote by $D'$ the boundary of
$U_N$, we suppose that the image of $D'$ by $p'$ is $D$. We obtain
the following proposition:

{\bf Proposition.} {\it The group $\hat{{Symp(N)}}_0$ acts naturally
on $V^1_{h_N}$. The stabilizer
 of the image of $\hat N\times N_U-D'$ by the covering map $p":\hat
N\times U_N\rightarrow V_{h_N}$ is $\hat{Symp(N,N_U)}_0$, the
connected component of $\hat{Symp(N)}_0$ which fixes $N_U$.}

\medskip

We suppose that the manifold $V$ is endowed with a symplectic form
$\omega_V$ such that the action of $H$ on $V$ is symplectic, and the
maps $\phi_i:U_i\rightarrow V$ which define the $(V,H)$-structure of
$N$ are symplectic maps. This enables to define on $V^1_{h_N}$ the
following symplectic structure: on the trivialization $U_i\times V$,
we define the form $\Omega_i$ to be the product of the restriction
of $\omega$ to $U_i$ with $-\omega_V$. Since the action of $H$ on
$V$ is symplectic, we deduce that the forms $\Omega_i$ glue together
to define a form $\Omega$ on $V^1_{h_N}$. Since we can assume that
the section $l_{h_N}$ takes it values in $V^1_{h_N}$, we have the
following proposition:

{\bf Proposition} {\it The section $l_{h_{N}}$ which defines the
$(V,H)$-structure is a Lagrangian submanifold. The orbits of
$l_{h_N}$ under $\hat{Symp(N)}_0$ are also Lagrangian submanifolds,
conversely every Lagrangian submanifold in the connected component
of $l_{h_N}$ in the space of Lagrangian submanifolds transverse to
the vertical foliation of $V^1_{h_N}$ is the image of $l_{h_N}$ by
an element of $\hat{Symp(N)}_0$.}
\medskip

{\bf proof.}

The section $l_{h_N}$ is defined on $U_i$ by
$l^i_{h_N}(u)=(u,\phi_i(u))$. Since the maps $\phi_i$ are symplectic
maps, we deduce that the section ${l^i}_{h_N}$ is a Lagrangian
submanifold as is the section $\phi(l^i_{h_N})$ defined by
$\phi(l^i_{h_N})(u)=(u,\phi\circ \phi_i(u))$.

Let $d'$ be a Lagrangian submanifolds of $V^1_{h_N}$ transverse to
the vertical foliation of $V^1_{h_N}$, $d'$ is the image of the
section $d$ of $V^1_{h_N}$. Consider a trivialization $(U_i\times
\phi_i(U_i),\phi_j\circ{\phi_i}^{-1})$ of $V^1_{h_N}$, we define the
symplectomorphism $\phi_d$ of $N$ whose restriction $\phi^i_d$ to
$U_i$ is defined for every element $u$ of $U_i$ by the projection of
$d_i(u)$ by the covering map $p":U_N\rightarrow N$, where $d_i$ is
the restriction of $d$ to $U_i$. The morphism $\phi_d$ is
well-defined indeed, suppose that $u\in U_i\cap U_j$, then
$d_j(u)=\phi_j\circ{\phi_i}^{-1}d_i(u)$. This implies that
$\phi^i_d(u)=\phi^j_d(u)$. The section $d$ is the image of $l_{h_N}$
by an element $\hat{\phi_d}$ of $\hat{Symp(N)}$ above $\phi_d$. If
$d'$ is in the connected component of $l_{h_N}$ in the space of
Lagrangian submanifolds transverse to the vertical foliation of
$V'_{h_N}$, we can suppose that $\phi_d$ is an element of
$Symp(N)_0$.
\medskip

Suppose that there exists a pseudo-complex structure $J$ adapted to
the symplectic form of $V^1_{h_N}$ such that the induced
differentiable metric is complete, $J$ enables to define an
isomorphism of bundle between the cotangent bundle of $l_{h_N}$,
$T^*l_{h_N}$ and $Nl_{h_N}$ the orthogonal of the tangent bundle
$Tl_{h_N}$  of $l_{h_N}$ in $V^1_{h_N}$,   We consider the
differentiable map $P_N:T^*l_{h_N}\rightarrow V^1_{h_N}$ defined as
follows: identify $T^*l_{h_N}$ with the normal bundle of $l_{h_N}$
in $V^1_{h_N}$ by the map which assigns to an element $n_u$ of the
fiber of the normal bundle at $u$ the $1$-form $i_{n_u}\Omega$, Let
$v_u$ be  an element of the fiber of $u$ in $T^*l_{h_N}$,
$P_N(v_u)=exp_u(v_u)$, where the exponential map is defined by the
differentiable metric $\Omega(J,.)$.

\bigskip

{\bf Theorem.}

{\it Suppose that $P_N$ is a symplectomorphism which induces a one
to one map between the Lagrangian submanifolds of
  $T^*l_{h_N}$ transverse to the vertical foliation of $T^*l_{h_N}$, and the Lagrangian submanifolds of
  $V^1_{h_N}$ transverse to the verical foliation of $V^1_{h_N}$, then the
orbit of $l_{h_N}$ under $\hat{Symp(N)}_0$ is contractible.}

\medskip

{\bf proof.}

Consider the cotangent bundle $T^*l_{h_N}$ $l_{h_N}$ endowed with
the differential of the Liouville form.  The map $P_N$ is a
symplectomorphism which send Lagrangian submanifolds of $T^*l_{h_N}$
transverse to the vertical foliation to Lagrangian submanifolds of
$V^1_{h_N}$ transverse to the vertical foliation. Since the
Lagrangian submanifolds of $T^*l_{h_N}$ transverse to the vertical
foliation are one to one with closed $1$-forms defined on $l_{h_N}$,
we deduce that the orbit of $l_{h_N}$ under $\hat{Symp(N)}_0$ is
contractible.

\bigskip

{\bf Symplectomorphisms group of symplectic affine manifolds.}

\medskip

Let $(N,\nabla_N)$ be a complete compact $n$-dimensional affine
manifold endowed with the parallel symplectic form $\omega$. The
bundle $V_{h_N}$ is the quotient of ${\R}^n\times{\R}^n$ by the
action of $\pi_1(N)$ which acts on the both factors by the holonomy
representation. The principal bundle associated to $V_{h_N}$ is the
bundle of affine frames.

{\bf Theorem} {\it Let $N$ be an $n$-dimensional compact affine
manifold, suppose that the affine structure of $N$ is defined by a
flat differentiable metric, and $N$ is endowed with a symplectic
form parallel in respect to the flat connection,  then if
$Ham(N,D)_0$ is  contractible, then $Ham(N)_0$ is also
contractible.}

{\bf proof.}

The manifold $N$ is the quotient of ${\R}^n$ by a subgroup
$\pi_1(N)$ whose linear part is contained in $O(n)$. The parallel
symplectic form $\omega_0$ is invariant by the holonomy and gives
rise to the form $\omega$ of $N$. The parallel complex structure
$J_0$ adapted to the flat metric of ${\R}^n$, and to $\omega_0$
gives rise to the pseudo-complex structure $J$ of $(N,\omega)$. In
this situation, $P_N:T^*l_{h_N}\rightarrow V_{h_N}$ is a
symplectomorphism.

Let $\alpha$ be a closed $1$-form defined on $N$, suppose that
$P_N(\alpha)$ is transverse to the vertical foliation, then it
defines an symplectomorphism $\hat{\phi}_{\alpha}$ of $V_{h_N}$
which gives rise to a symplectomorphism $\phi_{\alpha}$ of $N$, we
can define $\phi_{\alpha}(D)$  the image of $D$ by $\phi_{\alpha}$.
Suppose that $P_N(\alpha)$ is not transverse to the vertical
foliation, then the infinite dimensional Sard lemma implies the
existence of a sequence of closed $1$-forms $(\alpha_n)_{n\in{\N}}$,
such that $\alpha_n$ is transverse to the vertical foliation, and
$\alpha_n$ converges towards $\alpha$. We define $\phi_{\alpha}(D)$
to be the limit of $\phi_{\alpha_n}(D)$. This make sense because the
group of symplectomorphisms is $C^0$-closed in the group of
diffeomorphisms.

Let $N'(D)$ be the image of the space of closed $1$-forms
$C^1(\alpha)$ by the map $\Phi$ defined by
$\Phi(\alpha)=\phi_{\alpha}(D)$, $N'(D)$ contains an open and dense
subset of $N(D)$. Indeed, consider the subspace of image of closed
$1$-forms $\alpha$ such that $\alpha$ is transverse to the vertical
foliation. This image is open, since the neighborhood of a
symplectomorphism can be identified with a neighborhood of the zero
section in the cotangent bundle, and the image of the set $N"(D)$ of
closed $1$-forms such that the  image of their graph by $P_N$ are
transverse to the vertical foliation is dense in the space of graphs
of symplectomorphisms this implies also that  $N'(D)$ is closed in
$N(D)$, thus $N'(D)=N(D)$.

The space of closed $1$-forms ${C'}^1(N)$ which induces
symplectomorphisms which fix $D$ is a vector subspace since the
components of $D$ are affine submanifolds. To see this we can use a
generating functions-like theory. Indeed, let $\alpha$ be a closed
$1$-form, the image of the lift of the graph to the universal cover
of $\hat N$ by $P_N$ is the space of elements of ${\R}^n\times
{\R}^n$ of the form $(u_1+\alpha_1,..,u_n+\alpha_n
,u_1-\alpha_1,...,u_n-\alpha_n )$ where $(\alpha_1,..,\alpha_n)$ are
the component of the lift $\hat {\alpha}$ of $\alpha$ to ${\R}^n$.
The fact that $\phi_{\alpha}$ preserves $ D$ implies that if
$(u_1+\alpha_1 ,...,u_n+\alpha_n )$ is an element of the subspace
$\hat{D}$ of ${\R}^n$ over $ D$, then
$(u_1-\alpha_1,...,u_n-\alpha_n)$ is also an element of $\hat{ D}$.
 Let $\alpha$, $\alpha'$ be closed  $1$-forms
which define symplectomorphisms of $N$ which preserves $D$, then if
$U$ is a connected component of $\hat{ D}$ in ${\R}^n$, then the
element $(u_1-\alpha_1 ,...,u_n-\alpha_n )$ is tangent to $U$ if
$(u_1+\alpha_1,...,u_n+\alpha_n)$ in this situation, $(u_1,..,u_n)$
and $(\alpha_1(u),...,\alpha_n(u))$ are tangent to $U$. This implies
that $(\alpha_1+\alpha_1',...,\alpha_n+\alpha'_n)$ is tangent to
$U$, and henceforth that if $\alpha+\alpha'$ induces a
symplectomorphism  then the fact that
$(u_1-(\alpha_1+\alpha'_1),...,u_n-(\alpha_n+\alpha'_n))$ is an
element of $U$, if
$(u_1+\alpha_1+\alpha'_1,...,u_n+\alpha_n+\alpha'_n)$ is an element
of $U$ implies that this symplectomorphism preserves $D$. We have a
fibration ${C'}^1(N)\rightarrow C^1(N)\rightarrow N(D)$. Since
$C^1(N)$ and ${C'}^1(N)$ are contractible, we deduce that $N(D)$ is
contractible.

\medskip

{\bf Corollary.}

{\it The universal cover $\hat{Symp(N)}$ of the connected component
ot the group of symplectomorphisms $Symp(N)$ of a $4$-dimensional
compact manifold endowed with a flat metric, and a parallel
symplectic structure $\omega$   contractible.}

\medskip

{\bf Proof.} $N$ is the quotient of ${\R}^4$ by $\pi_1(N)$, the
boundary $D$ of this action is a configuration. We have a fibration
$Symp(N,D)_0\rightarrow Symp(N)_0\rightarrow N(D)$, since $N(D)$ is
contractible, we deduce that $Symp(N)_0$ has the homotopy type of
$Symp(N,D)_0$. We have an exact sequence $\hat{Ham(N,D)}\rightarrow
\hat{Symp(N,D)}\rightarrow L$, where $L$ is the image of
$\hat{Symp(N,D)}$ by the flux. Since $\hat{Ham(N,D)}$ is
contractible, we deduce that $\hat{Symp(N,D)}$ is  contractible.

\medskip

\centerline{\bf References.}

1. Abreu, McDuff, Topology of the symplectomorphisms group of
rational surfaces. J. of American Math Society.

2. Banyaga Sur la structure du groupe des diffeomorphismes qui
preservent une forme symplectiques. Comment. Math. Helvetici.

3. Joseph Coffey, Symplectomorphisms groups and isotropic skeletons,
Geometry and Topology.

4. Lalonde, Pinsonnault, The topology of the space of symplectic
balls in rational $4$-manifolds. Duke. Math. Journal.

5. Goldman, The geometric structures and varieties of
representations.

6. Pinsonnault, Remarques sur le groupe des symplectomorphismes de
l'eclatement de $S^2\times S^2$. These Universite U.Q.A.M.

7. Ruan, Y. Symplectic topology and extremal rays, Geometry and
functional analysis.

8. Tsemo, A. Homotopy of symplectomorphisms groups. In preparation.

\end{document}